\documentclass[12pt,thmsa]{article}
\usepackage{amsmath, latexsym, amsfonts, amssymb, amsthm, amscd}
\usepackage{graphics}
\newtheorem{theo}{Theorem}

\newtheorem{prop}{Proposition}
\newtheorem{lemma}{Lemma}

\newcommand{\ind}{\mbox{\rm 1\hspace{-0.04in}I}}

\begin{document}
\title{An $f$-divergence approach for optimal portfolios  in exponential Levy models}

\maketitle

\begin{center}
{\large S. Cawston}\footnote{$^{,2}$ LAREMA, D\'epartement de
Math\'ematiques, Universit\'e d'Angers, 2, Bd Lavoisier - 49045,
\\\hspace*{.4in}{\sc Angers Cedex 01.}

\hspace*{.05in}$^1$E-mail: suzanne.cawston@univ-angers.fr$\;\;\;$
$^2$E-mail: lioudmila.vostrikova@univ-angers.fr}{\large
 and  L. Vostrikova$^2$}
\end{center}
\vspace{0.2in}

\begin{abstract} We present a unified approach to get  explicit formulas for utility maximising strategies in Exponential Levy models. This approach is related to  $f$-divergence minimal martingale measures and based on a new concept of preservation of the Levy property by $f$-divergence minimal  martingale measures. For common $f$-divergences, i.e. functions which satisfy $f''(x)= ax^ {\gamma},\, a>0, \, \gamma \in \mathbb R$, we give the conditions for the existence of corresponding $u_f$- maximising strategies, as well as explicit formulas.
\end{abstract}

\section{Introduction}

Exponential Levy models have been widely used since the 1990's to represent asset prices. In the case of a.s. continuous trajectories, this leads to the classical Black-Scholes model, but the class of Levy models also contains a number of popular jump models including Generalized Hyperbolic models (\cite{E}) and Variance-Gamma models \cite{CGMY}. The use of such processes allows for an excellent fit both for daily log-returns (\cite{EK}) and intra-day data (\cite{EK}). The class is also flexible enough to allow for processes with either finite or infinite variation and finite or infinite activity. However, contrary to the Black-Scholes case, Levy models generally lead to incomplete financial markets : contingent claims cannot all be replicated by admissible strategies. Therefore, it is important to determine strategies which are, in a certain sense optimal. Various criteria are used, some of which are linked to risk minimisation ( see \cite{FS}, \cite{S1}, \cite{S2}) and others consisting in maximizing certain utility functions (see \cite{GR}, \cite{Ka}). 
It has been shown (see \cite{GR}, \cite{KSHA}) that such questions are strongly linked via the Fenchel-Legendre transform to dual optimisation problems on the set of equivalent martingale measures, i.e. the measures which are equivalent to the initial physical measure and under which the stock price is a martingale. More precisely, we recall that the convex conjugate of a concave function $u$ is defined by 
$$f(y)=\sup_{y\in\mathbb{R}}\{u(x)-xy\}=u(I(y))-yI(y)$$
where $I=(u')^{-1}$. In particular, we have the following correspondences : 
$$\begin{cases}
&\text{ if }u(x)=\ln(x) \text{ then }f(x)=-\ln(x)-1,
\\&\text{ if }u(x)=\frac{x^p}{p}, p<1, \text{ then }f(x)=-\frac{p-1}{p}x^{\frac{p}{p-1}},
\\&\text{ if }u(x)=1-e^{-x} \text{ then }f(x)=1-x+x\ln(x).
\end{cases}$$
Given a convex function $f$, the problem of minimising the $f$-divergence $E[f(\frac{dQ_T}{dP_T})]$ of the restrictions of the measures $P$ and $Q$ on the time interval $[0,T]$ over the set of equivalent martingale measures has been well studied for a number of functions  in \cite{ChS}, \cite{ChSL}, \cite{FM}, \cite{M},\cite{ES}, \cite{Kl}. For properties of $f$-divergence see also \cite{LV}. It has been noted in \cite{GR} that if a solution $Q^*$ to such a problem exists, there exists a predictable process $\hat{\phi}$ such that 
$$-f'(\frac{dQ^*_T}{dP_T})=x+\int_0^T \hat{\phi}_sdS_s,$$
where the process $S$ which represents the risky asset is a semi-martingale and $x$ is a constant. Moreover,
under some assumptions, $\hat{\phi}$ will then define a $u$-optimal strategy. However, it is in general far from easy to obtain an explicit expression for $\hat{\phi}$, although results exist for a certain number of special cases. These special cases concern what we will call common $f$-divergences, i.e. functions $f$ such that $f''(x) = a x^{\gamma}$ where $a>0$. \\
Our aim here is to obtain, for a certain class of utility functions, an explicit expression for $\hat{\phi}$ both when the Gaussian part of the Levy process is non-zero, i.e $c\neq 0$, and when $c=0$. We consider a class of $f$-divergences  whose $f$-divergence minimal martingale measure $Q^*$ preserves the Levy property of the initial Levy process. It is known that common $f$-divergences have preservation Levy property and the last class of $f$-divergences is larger then common $f$-divergences as it was shown in \cite{CV}. In addition, these new approach permit us to suggest a unified way for finding $\hat{\phi}$. In particular, we deduce from this result a unified formula for $\hat{\phi}$ for all common $f$-divergences.\\
Let us denote by $Z_T= \frac{dQ^*_T}{dP_T}$ the Radon-Nikodym derivative of $Q^*_T$ with respect to $P_T$ and let $(\beta , Y)$ be the Girsanov parameters for the change of measure from $P_T$ to $Q^*_T$ (cf. \cite{JSh}, p. 159). We consider utility functions $u$ such that their convex conjugate $f_u$ used as an $f$-divergence
gives us a Levy property preserving $f$-divergence minimal equivalent martingale measure $Q^*$. Then, under some integrability conditions, we  prove that
if the Gaussian part of the initial Levy process is not zero, then the optimal strategy $\hat{\phi}$ is given by:
$$\hat{\phi}^{(i)}_s=-\frac{\lambda\beta^{(i)} \,Z_{s-}}{S^{(i)}_{s-}}E_{Q^*}[f''(\lambda xZ_{T-s})Z_{T-s}]\,_{|_{x=Z_{s-}}}$$
where  $\lambda>0$ is the unique solution to the equation $E_{Q^*}[-f'(\lambda Z^*_T)]=x$ and $x$ is the initial capital.
If the Gaussian part of the initial Levy process is zero and the support of the Levy measure is of non-empty interior, then
$$ \hat{\phi}^{(i)} =-\frac{\lambda\gamma^{(i)} \,Z_{s-}}{S^{(i)}_{s-}}E_{Q^*}[f''(\lambda xZ_{T-s})Z_{T-s}]\,_{|_{x=Z_{s^-}}}$$
where $\gamma^{(i)}$ are   constants  related with the second Girsanov parameter and given by (\ref{gamma})(cf. Theorem 2).\\
In the particular case of common utility functions (corresponding to common $f$-divergences) we give conditions that ensure existence of the optimal strategy and its expression. For example, for $c\neq 0$,
$$\hat{\phi}^{(i)}_s=\frac{\alpha_{\gamma+1}(x)\,\beta^{(i)}}{E_{Q^*}[Z^{\gamma+1}_{s}]}\,\,\frac{\,Z^{\gamma+1}_{s-}}{S^{(i)}_{s-}}$$
where $\alpha_{\gamma+1}(x)$ is given by (\ref{alpha})(cf. Proposition 1).\\
The paper is organized in the following way: in 2. we recall known facts about utility maximisation, in 3. we prove (cf. Theorem 1) a decomposition needed to find optimal strategies, then in 4. we give a general result about optimal strategies, finally, in Proposition 1 we obtain the results concerning common $f$-divergences.
\section{Utility maximisation in exponential Levy models}
We start by describing our model in more detail. We assume that the financial market consists of a non-risky asset $B$ whose value at time $t$ is 
$$B_t=B_0e^{rt},$$
 where $r\geq 0$ is the interest rate which we assume to be constant, and $d$ risky assets whose prices are described by a $d$-dimensional stochastic process $S=(S_t)_{t\geq 0}$ with
$$S_t=(e^{X_t^{(1)}},...,e^{X_t^{(d)}})$$
where $X = (X_t^{(1)},...,X_t^{(d)})_{t\geq 0}$ is a $d$-dimensional Levy process defined on a filtered probability space $(\Omega,\mathcal{F},\mathbb{F},P)$ with the natural filtration $\mathbb F = (\mathcal F _t)_{t\geq 0}$ satisfying usual properties. We recall that Levy processes form the class of c\`{a}dl\`{a}g processes with stationary and independent increments and such that the law of $X_t$ is given by the Levy-Khintchine formula : for all $t\geq 0$, for all $u\in\mathbb{R}$
 $$E[e^{i<u,X_t>}]=e^{t\psi(u)}$$
with 
$$\psi(u)=i<u,b>-\frac{1}{2}\,^{\top}ucu+\int_{\mathbb{R}^d}[e^{i<u,y>}-1-i<u,h(y)> ]\nu(dy)$$
where $b\in\mathbb{R}^d$ is a drift, $c$ is a positive $d\times d$ symmetric matrix, $\nu$ is a positive measure on $\mathbb{R}^{d}\setminus \{0\}$ which satisfies 
$$\int_{\mathbb{R}^d}1\wedge |y|^2 \nu(dy)<+\infty $$
and $h(\cdot )$ is a truncation function.
The triplet $(b,c,\nu)$ entirely determines the law of the Levy process $X$, and is called the characteristic triplet of $X$. For more details see  \cite{Sa}. We also recall that if $S=e^X$, there exists a Levy process $\hat{X}$ such that $S=\mathcal{E}(\hat{X})$, where $\mathcal{E}$ denotes the Doleans-Dade exponential. For more details see \cite{JSh}.\\

An investor will share out his capital among the different assets according to a strategy which is represented by a process $\Phi=(\eta,\phi)$, where $\eta$ represents the quantity invested in the non-risky asset $B$, and $\phi=(\phi^{(1)},...,\phi^{(d)})$ is the quantity invested in the risky assets. From now on, we will denote by 
$$(\phi\cdot S)_t= \sum_{i=1}^d \int_0^t \phi_s^{(i)} dS_s^{(i)}$$
the variation of capital due to the investment in the risky assets. We now define more precisely our set of admissible strategies. 
We recall that an admissible strategy is a predictable process $\Phi=(\eta,\phi)$ taking values in $\mathbb{R}^{d+1}$, such that $\eta$ is $B$-integrable, $\phi$ is $S$-integrable and for which there exists $a\in\mathbb{R^+}$ such that for all $t\geq 0$, 
$$(\phi\cdot S)_t\geq -a$$
We denote by $\mathcal{A}$ the set of all admissible strategies.\\ 

We are interested in strategies which are optimal in the sense of utility maximisation. We recall that a utility function is a function $u:]\underline{x},+\infty[\longrightarrow \mathbb{R}$, which is $\mathcal{C}^1$, strictly increasing, strictly concave and such that 
$$\lim_{x\to+\infty}u'(x)=0 \text{ and }\lim_{x\to\underline{x}}u'(x)=+\infty$$
where $\underline{x}=\inf\{ x \,|\,x \in dom(u)\}$. In particular, the most common utility functions are $u(x)=\ln(x)$, $u(x)=\frac{x^p}{p}$, $p<1$, or $u(x)=1-e^{-x}$. 
We now recall the definition of u-optimal and u-asymptotically optimal strategies. This last notion was first introduced in \cite{Ka}.  It will allow  us to consider in a  unified way all utilities including those with $\underline{x}=\-\infty$. \\

\noindent We say that a strategy $\hat{\phi}\in\mathcal{A}$ is $u$-optimal on $[0,T]$ if 
$$E[u(x+(\hat{\phi}\cdot S)_T)]=\sup_{\phi\in\mathcal{A}}E[u(x+(\phi\cdot S)_T)]$$
A sequence of admissible strategies $(\hat{\phi}^{(n)})_{n\geq 1}$ is asymptotically $u$-optimal on $[0,T]$ if 
$$\lim_{n\to+\infty}E[u(x+(\hat{\phi}^{(n)}\cdot S)_T)]=\sup_{\phi\in\mathcal{A}}E[u(x+(\phi\cdot S)_T)].$$

\section{A decomposition for Levy preserving  equivalent martingale measures}
In this section, we consider a fixed strictly convex function $f$, $f\in C^3(\mathbb R^{+,*})$, and a Levy preserving equivalent  martingale measure $Q$ whose density is given by the process $Z= (Z_t)_{t\geq 0}$. We recall that $Q$ preserves the Levy property if
$X$ remains a Levy process under $Q$. We also recall that we characterize the change of measure from $P$ into $Q$ by the Girsanov parameters $(\beta, Y)$. Then the fact that $Q$ preserves the Levy property can be seen as a change of measure such that the first Girsanov parameter $\beta$ is a constant and the second parameter $Y$  depends only on jump-sizes. As a consequence, the density of a Levy preserving measure is of the form $Z=\mathcal{E}(N)$, where 
$$N_t=\beta X^{(c)}_t+\int_0^t \int_{\mathbb{R}^{d*}}(Y(x)-1)(\mu^{X}-\nu ^{X,P})(ds,dx)$$
In addition, if $Q$ is a martingale measure then
 $\beta$ and $Y$ satisfy
$$b+\frac{1}{2} diag(c)+c\beta+\int_{\mathbb{R}^d}[(e^x-1)Y(x)-h(x)]\nu(dx)=0$$
The last relation  ensures that the drift of $S$ under the measure $Q$ is zero.\\
\par Our main aim in this section is to show that under certain integrability conditions, the decomposition given in Theorem \ref{thdec} holds.
 We introduce c\`{a}dl\`{a}g versions of the processes $(\xi_t(x))_{t\geq 0}$ et $(H_t(x, y))_{t\geq 0}$
where for $0\leq t \leq T$
\begin{equation}\label{xi}
\xi_t(x)= E_Q[f''(xZ_{T-t})Z_{T-t}]
\end{equation}
and 
\begin{equation}\label{H}
 H_t(x,y)=E_Q[f'(xZ_{T-t}Y(y))-f'(xZ_{T-t})]
\end{equation}
\begin{theo}\label{thdec}
 Let $f$ be a strictly convex function belonging to $ C^3(\mathbb R^{+,*})$. 
Let $Z$ be the density of a Levy preserving equivalent martingale measure $Q$. Assume that 
$Q$ satisfies : for all $\lambda >0$ and all compact set $K\subseteq \mathbb R^+$
\begin{equation}\label{integcd}
E_P|f(\lambda Z_T)|<+\infty,\,\,\,\,E_Q|f'(\lambda\,Z_T)|<+\infty,\,\,\,\,\sup_{t\leq T}\sup_{\lambda \in K}E_Q[f''(\lambda Z_t)Z_t]<+\infty.
\end{equation}
 Then, for all $\lambda>0$ we have $Q$- a.s, for all $t\leq T$,  
\begin{equation}\label{decf}
E_Q[f'(\lambda Z_T)|\mathcal{F}_t]=E_Q[f'(\lambda Z_T)] +
\end{equation}
$$\sum_{i=1}^d \lambda \beta ^{(i)}\int_0^t \xi_s(\lambda Z_{s-})\, Z_{s-}\,d X_s ^{(c),Q,i}+\int_0^t \int_{\mathbb{R}^d}H_s(\lambda Z_{s-},y)\,(\mu ^X - \nu ^{X, Q})(ds, dy)$$
\end{theo}
This result is based on an application of the Ito formula, but it will require some technical lemmas.

\vskip 0.2cm We recall that as $Q$ preserves the Levy property, for all $t\leq T$, $Z_t$ and $\frac{Z_T}{Z_t}$ are independent under $P$ and that $\mathcal{L}(\frac{Z_T}{Z_t}\,|\,P)=\mathcal{L}(Z_{T-t}\,|\,P)$. Therefore 
$$E_Q[f'(\lambda Z_T)|\mathcal{F}_t]=\rho(t,Z_t)$$
where $\rho(t,x)=E_Q[f'(\lambda xZ_{T-t})]$. Our integrability conditions do not allow us to apply the Ito formula directly to the function $\rho(t,Z_t)$. Therefore, we start by considering a sequence of bounded approximations of $f'$, and will then obtain (\ref{decf}) by studying the convergence of
analogous decompositions for the approximations of $f'$.
 
\begin{lemma}\label{approx} Let $f$ be a strictly convex function belonging to $ C^3(\mathbb R^{+,*})$. There exists a sequence of bounded increasing functions $(\phi_n)_{n\geq 1}$, which are of class $\mathcal{C}^2$ on $\mathbb{R}^{+,*}$, such that for all $n\geq 1$, $\phi_n$ coincides with $f'$ on the compact set $[\frac{1}{n},n]$ and such that for $n$ large enough and  for all $x,y>0$   the following inequalities hold~:
\begin{equation}\label{apprfprime}
|\phi_n(x)|\leq 4|f'(x)|+ \alpha \text{ , }|\phi'_n(x)|\leq 3f''(x) \text{ , }|\phi_n(x)-\phi_n(y)|\leq 5|f'(x)-f'(y)|
\end{equation}
where $\alpha$ is a real positive constant.
\end{lemma} 

\noindent \it Proof \rm 
We set, for $n\geq 1$, 
$$A_n(x)=f'(\frac{1}{n})-\int_{x\vee \frac{1}{2n} }^{\frac{1}{n}}f''(y)(2ny-1)^2(5-4ny)dy$$
$$B_n(x)=f'(n)+\int_n^{x\wedge (n+1)} f''(y)(n+1-y)^2(1+2y-2n)dy$$
and finally
$$\phi_n(x)=\begin{cases}
& A_n(x) \text{ if }0 \leq x< \frac{1}{n},
\\&f'(x) \text{ if }\frac{1}{n}\leq x\leq n,
\\& B_n(x) \text{ if }x>n.

\end{cases}$$

Here $A_n$ and $B_n$ are defined so that $\phi_n$ is of class $\mathcal{C}^2$ on $\mathbb{R}^{+,*}$.
For the inequalities we use the fact that $f'$ is increasing function and the estimations:
$0\leq (2nx-1)^2(5-4nx)\leq 1$ for $\frac{1}{2n} \leq x\leq \frac{1}{n}$ and
$0\leq(n+1-x)^2(1+2x-2n)\leq 3$ for $n \leq x\leq n+1$. 
$\Box$

\vskip 0.2cm We now introduce for each $n\geq 1$ the function 
$$\rho_n(t,x)=E_Q[\phi_n(\lambda xZ_{T-t})].$$
and we obtain the following analog to Theorem \ref{thdec}, replacing  $f'$ with $\phi_n$.
For that we introduce for $0\leq t\leq T$
\begin{equation}\label{xin}
\xi^{(n)}_t(x)= E_Q[\phi_n'(xZ_{T-t})Z_{T-t}]
\end{equation}
and
\begin{equation}\label{han}
 H^{(n)}_t(x,y)=E_Q[\phi_n(xZ_{T-t}Y(y))-\phi_n(xZ_{T-t})]
\end{equation}
 
\begin{lemma}\label{decomp}
We have $Q$-a.s., for all $t\leq T$,
\begin{equation}\label{ndecomp}
\rho_n(t,Z_t)=E_Q[\phi_n(\lambda Z_T)] +
\end{equation}
$$\sum_{i=1}^d  \lambda\,\beta ^{(i)} \int_0^t  \xi^{(n)}_s (\lambda Z_{s-})\,Z_{s-}\,dX^{(c),Q,i}_s+\int_0^t \int_{\mathbb{R}^{d}}H^{(n)}_s(\lambda Z_{s-},y)\,(\mu^X-\nu ^{X,Q})(ds, dy)$$
where $\beta = \,^\top\!(\beta _1,\cdots, \beta _d) $ and $\nu ^{X,Q}$ is the dual predictable projection or the compensator of the jump measure $\mu^X$ with respect to $(\mathbb F, Q)$.
\end{lemma}
\noindent \it Proof \rm In order to apply the Ito formula to $\rho_n$, we need to show that $\rho_n$ is twice continuously differentiable with respect to $x$ and once with respect to $t$ and that the corresponding derivatives are bounded for all $t\in[0,T]$ and $x\geq \epsilon$, $\epsilon >0.$ First of all, we note from the definition of $\phi_n$ that 
for all $x\geq \epsilon >0$
$$|\frac{\partial}{\partial x}\phi_n(\lambda xZ_{T-t})| = |\lambda Z_{T-t}\phi'_n(\lambda xZ_{T-t})|\leq \frac{(n+1)}{\epsilon}\sup_{z>0}|\phi'_n(z)|<+\infty.$$
Therefore, $\rho_n$ is differentiable with respect to $x$ and we have
$$\frac{\partial}{\partial x}\rho_n(t,x)=\lambda E_Q[\phi'_n(\lambda xZ_{T-t})\,Z_{T-t}].$$
Moreover, the function $(x,t)\mapsto \lambda\phi'_n(\lambda xZ_{T-t})Z_{T-t}$ is continuous $P$-a.s. and bounded. This implies that $\frac{\partial}{\partial x}\rho_n$ is continuous and bounded for
$t\in [0,T] $ and $x\geq \epsilon >0$. \\
 In the same way, for all $x\geq \epsilon >0$
$$|\frac{\partial^2}{\partial x^2}\phi_n(\lambda xZ_{T-t})|=\lambda^2 Z_{T-t}^2\phi''_n(\lambda x Z_{T-t})\leq \frac{(n+1)^2}{\epsilon ^ 2} \sup_{z>0}\phi''_n(z)<+\infty .$$
Therefore, $\rho_n$ is twice continuously differentiable in $x$ and 
$$\frac{\partial^2}{\partial x^2}\rho_n(t,x)=\lambda^2 E_Q[\phi_n''(\lambda xZ_{T-t})Z^2_{T-t}]$$
We can verify easily that it is again continuous and bounded function.
In order to obtain differentiability with respect to $t$, we need to apply the Ito formula to $\phi_n$ :
$$\begin{aligned}
\phi_n(\lambda xZ_{t})=&\phi_n(\lambda x)+\sum_{i=1}^d \int_0^{t}\lambda x\phi'_n(\lambda xZ_{s-})\beta ^{(i)}\,Z_{s-}dX^{(c),Q,i}_s
\\&+\int_0^{t}\int_{\mathbb{R}^{d}}[\phi_n(\lambda xZ_{s-}Y(y))-\phi_n(\lambda xZ_{s-})]\,(\mu^X-\nu^{X,Q})(ds, dy)
\\&+\int_0^{t}\psi_n(\lambda x,Z_{s-})ds
\end{aligned}$$
where 
$$\begin{aligned}
\psi_n(\lambda x,Z_{s-})=& ^{\top}\beta c\beta [\lambda x Z_{s-}\phi'_n(\lambda xZ_{s-})+\frac{1}{2}\,x^2\lambda ^ 2 Z_{s-}^ 2 \phi''_n(\lambda xZ_{s-})]
\\&+\int_{\mathbb{R}^{d}}[(\phi_n( \lambda xZ_{s-}Y(y))-\phi_n(\lambda xZ_{s-}))\,Y(y)-\lambda x\phi'_n(\lambda xZ_{s-})Z_{s-}(Y(y)-1)]\nu(dy).
\end{aligned}$$
Therefore, for fixed $t>0$ 
$$E_Q[\phi_n(\lambda xZ_{T-t})]=\int_0^{T-t}E_Q[\psi_n(\lambda x,Z_{s-})]ds$$
so that $\rho_n$ is differentiable with respect to $t$ and 
$$\frac{\partial}{\partial t}\rho_n(t,x)=-E_Q[\psi_n(\lambda x,Z_{s-})]_{|s=T-t}.$$
We can also easily check that this is again a continuous and bounded function. For this we use the fact that $ \phi_n$, $\phi'_n$ and  $\phi''_n $  are bounded functions and also that the Hellinger process
of $Q_T$ and $P_T $ of the order $1/2$  is finite.
\\ 

We can finally apply the Ito formula to $\rho_n$. For that we use the stopping times
$$s_m= \inf \{ t\geq 0\,|\, Z_t \leq \frac{1}{m} \},$$
with $m\geq 1$ and $\inf\{\emptyset\}= +\infty$.
Then, from Markov property of Levy process we have :
$$\rho _n(t\wedge s_m, Z_{t\wedge s_m}) = E_Q(\phi _n(\lambda Z_T)\,|\, \mathcal F_{t\wedge s_m})$$
We remark that $(E_Q(\phi _n(\lambda Z_T)\,|\, \mathcal F_{t\wedge s_m})_{t\geq 0}$ is $Q$-martingale, uniformly integrable with respect to $m$. From Ito formula we have :
\begin{eqnarray*} \rho _n(t\wedge s_m, Z_{t\wedge s_m}) = E_Q(\phi _n(\lambda Z_T)) + \int_0^{t\wedge s_m}\frac{\partial \rho _n}{\partial s}(s, Z_{s-}) ds +\\  \int_0^{t\wedge s_m}\frac{\partial \rho _n}{\partial x}(s, Z_{s-}) dZ_s +
 \frac{1}{2}\int_0^{t\wedge s_m}\frac{\partial^2 \rho _n}{\partial x^2}(s, Z_{s-}) d< Z^c>_s +\\
\sum _{0\geq s\geq t\wedge s_m}\rho _n(s, Z_{s})- \rho _n(s, Z_{s-}) - \frac{\partial \rho _n}{\partial x}(s, Z_{s-}) \Delta Z_s
\end{eqnarray*}
where $\Delta Z_s = Z_s - Z_{s-}$.
After some standard simplifications, we see that
$$\rho _n(t\wedge s_m, Z_{t\wedge s_m}) =   A_{t\wedge s_m} + M_{t\wedge s_m}$$
where $(A_{t\wedge s_m})_{0\leq t\leq T}$ is predictable process, which is equal to zero,
\begin{eqnarray*}A_{t\wedge s_m} = \int_0^{t\wedge s_m}\frac{\partial \rho _n}{\partial s}(s, Z_{s-}) ds +
 \frac{1}{2}\int_0^{t\wedge s_m}\frac{\partial^2 \rho _n}{\partial x^2}(s, Z_{s-}) d< Z^c>_s +\\
\int_0^{t\wedge s_m}\int _{\mathbb R}[ \rho _n(s, Z_{s-}+x)- \rho _n(s, Z_{s-}) - \frac{\partial \rho _n}{\partial x}(s, Z_{s-})x] \nu^ {Z,Q}(ds, dx)
\end{eqnarray*}
and $(M_{t\wedge s_m})_{0\leq t\leq T}$ is a $Q$-martingale,
\begin{eqnarray*}M_{t\wedge s_m} =  E_Q(\phi _n(\lambda Z_T)) +  \int_0^{t\wedge s_m}\frac{\partial \rho _n}{\partial x}(s, Z_{s-}) dZ^c_s+\\
\int_0^{t\wedge s_m}\int _{\mathbb R}[ \rho _n(s, Z_{s-}+x)- \rho _n(s, Z_{s-})](\mu^ Z(ds, dx)-\nu^ {Z,Q}(ds, dx))
\end{eqnarray*}
Then, we pass to the limit as $m\rightarrow +\infty$. We remark that the sequence $(s_m)_{m\geq 1}$ is going to $+\infty$ as $m\rightarrow \infty$. From \cite{RY}, corollary 2.4, p.59, we obtain that
$$\lim_{m\rightarrow \infty} E_Q(\phi_n(Z_T)\,|\, \mathcal F_{t\wedge s_m})= E_Q(\phi_n(Z_T)\,|\, \mathcal F_{t})$$ and by the definition of local martingales we get:
$$\lim_{m\rightarrow \infty} \int_0^{t\wedge s_m}\frac{\partial \rho _n}{\partial x}(s, Z_{s-}) dZ^c_s =  \int_0^{t}\frac{\partial \rho _n}{\partial x}(s, Z_{s-}) dZ^c_s = \int_0^{t}\lambda \xi_s^{(n)}(Z_{s-}) dZ^c_s$$
and
$$\lim_{m\rightarrow \infty}\int_0^{t\wedge s_m}\int _{\mathbb R}[ \rho _n(s, Z_{s-}+x)- \rho _n(s, Z_{s-})](\mu^ Z(ds, dx)-\nu^ {Z,Q}(ds, dx)) =$$
$$ \int_0^{t}\int _{\mathbb R}[ \rho _n(s, Z_{s-}+x)-\rho _n(s, Z_{s-})](\mu^ Z(ds, dx)-\nu^ {Z,Q}(ds, dx))$$
Now, in each stochastic integral we pass from the integration with respect to the process $Z$ to the one with respect to the process $X$. For that we remark that
$$dZ^c_s = \sum _{i=1}^d \beta ^{(i)} Z_{s-} dX_s^{c,Q,i},\,\,\,\Delta Z_s = Z_{s-}\Delta X_s,\,\,\, Y(\Delta X_s)=
 1+ \frac{\Delta X_s}{Z_{s-}}.$$ Lemma \ref{decomp} is proved.
$\Box $

\vskip 0.2cm
We now turn to the proof of Theorem \ref{thdec}. In order to obtain the decomposition for $f'$, we prove convergence in probability of the different processes which appear in (\ref{ndecomp}). 
\par \noindent \it Proof of Theorem \ref{thdec} \rm 
For   $n\geq 1$ and a fixed $\lambda >0$, we introduce the stopping times 
\begin{equation}\label{taun}
\tau_n= \inf\{t\geq 0\,|\,\lambda Z_t\geq n \text{ or }\lambda Z_t\leq \frac{1}{n}\}
\end{equation}
where $\inf\{\emptyset\}= +\infty$ and we note that $\tau _n\rightarrow +\infty$ ($P$-a.s.) as $n\rightarrow \infty\,$.
First of all, we note that 
$$|E_Q[f'(\lambda Z_T)|\mathcal{F}_t]-\rho_n(t,Z_t)|\leq E_Q[|f'(\lambda Z_T)-\phi_n(\lambda Z_T)||\mathcal{F}_t]$$
As $f'$ and $\phi_n$ coincide on the interval $[\frac{1}{n},n]$, it follows from Lemma \ref{decomp} that 
$$\begin{aligned}
|E_Q[f'(\lambda Z_T)|\mathcal{F}_t]-\rho_n(t,Z_t)|&\leq E_Q[|f'(\lambda Z_T)-\phi_n(\lambda Z_T)|{\bf 1}_{\{\tau_n<T\}}|\mathcal{F}_t]
\\&\leq  E_Q[(5|f'(\lambda Z_T)|+ \alpha){\bf 1}_{\{\tau_n<T\}}|\mathcal{F}_t].
\end{aligned}$$
Now, for every $\epsilon >0$, by the Doob's inequality and the Lebesgue dominated convergence theorem we get: 
$$\lim_{n\to+\infty}Q(\sup_{t\leq T}E_Q[(5|f'(\lambda Z_T)|+\alpha){\bf 1}_{\{\tau_n<T\}}|\mathcal{F}_t] > \epsilon )\leq \lim_{n\to+\infty}\frac{1}{\epsilon }E_Q[(5|f'(\lambda Z_T)|+\alpha){\bf 1}_{\{\tau_n<T\}}]=0$$
Therefore, we have 
$$\lim_{n\to+\infty}Q(\sup_{t\leq T}|E_Q[f'(\lambda Z_T) - \rho_n(t,\lambda Z_t)|\mathcal{F}_t]|>\epsilon)=0.$$
We now turn to the convergence of the three elements on the right-hand side of (\ref{ndecomp}). 
We have almost surely $\lim_{n\to+\infty}\phi_n(\lambda Z_T)=f'(\lambda Z_T)$, and for all $n\geq 1$, $|\phi_n(\lambda Z_T)|\leq 4|f'(\lambda Z_T)|+\alpha $. Therefore, it follows from the dominated convergence theorem that, 
$$\lim_{n\to+\infty}E_Q[\phi_n(\lambda Z_T)]=E_Q[f'(\lambda Z_T)].$$
 We now prove the convergence of the continuous martingale parts of (\ref{ndecomp}). It follows from Lemma \ref{approx} that 
$$\begin{aligned}                       
|\xi^{(n)}_t(\lambda Z_t)-\xi_t(\lambda Z_t)|\leq &  E_Q[Z_T|\phi'_n(\lambda Z_T)-f''(\lambda Z_T)|\mathcal{F}_t]\leq
\\& 4 E_Q[Z_T|f''(\lambda Z_T)|{\bf 1}_{\{\tau_n<T\}}|\mathcal{F}_t].
\end{aligned}$$                     
Hence, we have as before  for $\epsilon >0 $           
$$\lim_{n\to+\infty}Q(\sup_{t\leq T}|\xi^{(n)}_t(\lambda Z_t)-\xi_t(\lambda Z_t)|>\epsilon)\leq \lim_{n\to+\infty}\frac{4 }{\epsilon}E_Q[Z_Tf''(\lambda Z_T){\bf 1}_{\{\tau_n<T\}}]=0$$

 Therefore, it follows from the Lebesgue dominated convergence theorem for stochastic integrals (see \cite{JSh}, Theorem I.4.31, p.46 ) that for all $\epsilon >0$ and $1\leq i\leq d$ 
$$\lim_{n\to+\infty}Q(\sup_{t\leq T}\,\big|\int_0^t (\xi^{(n)}_s(\lambda Z_{s-})-\xi_s(\lambda Z_{s-}))dX^{(c),Q,i}_s \big|>\epsilon)=0.$$
It remains to show the convergence of the discontinuous martingales to zero as $n\rightarrow \infty $.
 We start by writing 
$$\int_0^{t}\int_{\mathbb{R}^{d}}[H^{(n)}_s(\lambda Z_{s-},y)-H_s(\lambda Z_{s-},y)](\mu^X-\nu^{X, Q})(ds, dy) =M^{(n)}_t+N^{(n)}_t$$
with 
$$M^{(n)}_t=\int_0^t \int_{ A}[H^{(n)}_s(\lambda Z_{s-},y)-H_s(\lambda Z_{s-},y)] (\mu^X-\nu^{X, Q})(ds, dy),$$
$$N^{(n)}_t=\int_0^t \int_{ A^ c}[H^{(n)}_s(\lambda Z_{s-},y)-H_s(\lambda Z_{s-},y)] (\mu^X-\nu^{X, Q})(ds, dy),$$
where $ A = \{y : |Y(y)-1|<\frac{1}{4}\}$.\\
For $p\geq 1$, we consider the sequence of stopping times $\tau_p$ defined by (\ref{taun})   with replacing $n$ by a real positive $p$. We also introduce the  processes $$M^{(n,p)}= (M^{(n,p)}_t)_{t\geq 0},\,\,
N^{(n,p)}= (N^{(n,p)}_t)_{t\geq 0}$$ with $M^{(n,p)}_t= M^{(n)}_{t\wedge \tau _p}$, $N^{(n,p)}_t= N^{(n)}_{t\wedge \tau _p}$.
We remark that for $p\geq 1$ and $\epsilon >0$
$$ Q(\sup _{t\leq T} |M^{(n)}_t+N^{(n)}_t| >\epsilon) \leq Q(\tau _p <T) + Q( \sup _{t\leq T} |M^{(n,p)}_t| >\frac{\epsilon }{2}) +
 Q( \sup _{t\leq T} |N^{(n,p)}_t| >\frac{\epsilon }{2}).$$
Furthermore, we obtain  from the Doob's martingale inequalities that
\begin{equation}\label{02}
 Q( \sup _{t\leq T} |M^{(n,p)}_t|  >\frac{\epsilon}{2}) \leq \frac{4}{\epsilon ^ 2}\mathbb E_Q [(M^{(n,p)}_T)^ 2]
\end{equation}
and
\begin{equation}\label{03}
 Q( \sup _{t\leq T} |N^{(n,p)}_t|  >\frac{\epsilon }{2}) \leq \frac{2}{\epsilon }\mathbb E_Q |N^{(n,p)}_T|
\end{equation} 
Since $\tau _p \rightarrow + \infty$ as $p\rightarrow +\infty$ it is sufficient to show  that $E_Q[M^{(n,p)}]^2$ and $E_Q|N^{(n,p)}|$ converge to $0$ as $n\rightarrow \infty $.\\
To do so we estimate $E_Q [(M^{(n,p)}_T)^ 2]$ and prove that\\\\ 
$E_Q [(M^{(n,p)}_T)^ 2]\leq$\\
$$C \big(\int _0^ T\sup _{v\in K}\, E_Q^2[ Z_{s}\,f''( v Z_{s})\ind_{\{\tau _{q_n}<s\}}]ds\big)\,\big(\int _ { A} (\sqrt{Y(y)}-1)^ 2\nu(dy) \big)$$
where $C$ is a positive constant, $K$ is some compact set of $\mathbb R^{+,*}$ and $q_n=\frac{ n }{4p }$.\\\\
First we note  that on stochastic interval $[\![0, T\wedge \tau _p)]\!]$ we have $1/p\leq \lambda Z_{s-}\leq p$, and, hence,
$$E_Q[(M_T^{(n,p)})^2]=E_Q[\int_0^{T\wedge \tau _p}\int_{ A}|H^{(n)}_s(\lambda Z_{s-},y)-H_s(\lambda Z_{s-},y)|^2\,Y(y)\nu(dy)ds]\leq$$
$$\int_0^T\int_{ A}\sup_{1/p\leq x\leq p}|H^{(n)}_{T-s}(x,y)-H_{T-s}(x,y)|^2\,Y(y)\nu(dy)ds]$$
To estimate the difference $|H^{(n)}_{T-s}(x,y)-H_{T-s}(x,y)|$ we note that
$$H^{(n)}_{T-s}(x,y)-H_{T-s}(x,y)= E_Q[\phi_n(xZ_{s}Y(y))-\phi_n(xZ_{s})-f'( xZ_{s}Y(y))+f'( xZ_{s})]$$
From  Lemma \ref{approx} we deduce that if $ x Z_{s} Y(y) \in [1/n, n]$ and $ x Z_{s} \in [1/n, n]$ then the expression on the right-hand side of the previous equality is zero. But if $y\in  A$ we also have : $1/4 \leq Y(y)\leq 5/4$ and, hence,\\
$|H^{(n)}_{T-s}(x,y)-H_{T-s}(x,y)|\leq$ \\ 
$$|E_Q[\ind_{\{\tau _{q_n}<s\}}|\phi_n(xZ_{s}Y(y))-\phi_n(xZ_{s})-f'(xZ_{s}Y(y))+f'(xZ_{s})|]$$
Again from the inequalities of Lemma \ref{approx} we get:
$$|H^{(n)}_{T-s}( x,y)-H_{T-s}( x,y)|\leq  6 E_Q[\ind_{\{\tau _{q_n}<s\}}|f'( xZ_{s}Y(y))-f'( xZ_{s})|]$$
Writing 
$$f'( xZ_{s}Y(y))-f'( xZ_{s})= \int _1^ {Y(y)} x Z_{s} f''( x Z_{s}\theta) d\theta$$
we finally get
$$|H^{(n)}_{T-s}(x,y)-H_{T-s}(x,y)|\leq  6 \,x\,|Y(y)-1|\,\sup_{1/4\leq u \leq 5/4}E_Q[\ind_{\{\tau _{q_n}<s\}}\, Z_{s}\,f''( x uZ_{s})]$$
and this gives us the estimation of $E_Q [(M^{(n,p)}_T)^ 2]$ cited above.\\
We know that $P_T\sim Q_T$ and this means that the corresponding Hellinger process of order 1/2 is finite:
$$ h_T(P, Q, \frac{1}{2}) =\frac{T}{2}\,^{\top}\beta c \beta  + \frac{T}{8} \int_{\mathbb R}(\sqrt{Y(y)}-1)^ 2\nu(dy) < +\infty .$$
Then 
$$\int_{ A}(\sqrt{Y(y)}-1)^ 2\nu(dy) < +\infty .$$
From the Lebesgue dominated convergence theorem  and (\ref{integcd}) we get:
$$\int _0^ T\sup _{v\in K}\, E_Q^2[ Z_{s}f''( v Z_{s})\ind_{\{\tau _{q_n}<s\}}]ds \rightarrow 0$$
as $n\rightarrow +\infty$ and this information together with the estimation of $E_Q[(M_T^{(n,p)})^2]$ proves the convergence of $E_Q[(M_T^{(n,p)})^2]$  to zero as $n\rightarrow +\infty$.

\par We now turn to the convergence of $E_Q| N_T^{(n,p)}|$ to zero as $n\rightarrow +\infty$. 
For this we prove that
$$E_Q| N_T^{(n,p)}|\leq 2 T E_Q[\ind_{\{\tau _{n}<T\}}(5|f'(\lambda Z_{T})|+\alpha)]\int_{ A^ c}Y(y) d\nu $$
We start by noticing that
$$E_Q|N^{(n,p)}_T |\leq 2 E_Q[\int_0^{T\wedge \tau _p}\int_{ A^c}|H^{(n)}_s(\lambda Z_{s-},y)-H_s(\lambda Z_{s-},y)|\,Y(y)\nu(dy)ds]\leq $$
$$ 2 \int_0^{T}\int_{ A^c}E_Q[|H^{(n)}_s(\lambda Z_{s-},y)-H_s(\lambda Z_{s-},y)|\,Y(y)\nu(dy)ds]$$
 To evaluate the right-hand side of the previous inequality we write\\
$|H^{(n)}_s(\lambda x,y)-H_s(\lambda x,y)|\leq$\\
 $$E_Q|\phi_n(\lambda xZ_{T-s}Y(y))-f'(\lambda xZ_{T-s}Y(y))|  + E_Q |\phi_n(\lambda xZ_{T-s})-f'(\lambda xZ_{T-s})|.$$
We remark that in law with respect to $Q$
$$ |\phi_n(\lambda x Z_{T-s}Y(y))-f'(\lambda x Z_{T-s}Y(y))| = E_Q[|\phi_n(\lambda Z_{T})-f'(\lambda Z_{T})|\,|\,Z_{s}=x\,Y(y)]$$
and
$$ |\phi_n(\lambda xZ_{T-s})-f'(\lambda xZ_{T-s})| = E_Q[|\phi_n(\lambda Z_{T})-f'(\lambda Z_{T})|\,|\,Z_{s}=x]$$
Then
$$|H^{(n)}_s(\lambda x,y)-H_s(\lambda x,y)| \leq 2 E_Q|\phi_n(\lambda Z_{T})-f'(\lambda Z_{T})|$$
From  Lemma \ref{approx} we get:
$$E_Q|\phi_n(\lambda xZ_{T})-f'(\lambda xZ_{T})|\leq E_Q[\ind_{\{\tau _{n}<T\}}|\phi_n(\lambda Z_{T})-f'(\lambda Z_{T})|]\leq E_Q[\ind_{\{\tau _{n}<T\}}(5|f'(\lambda Z_{T})|+\alpha)]$$
and is proves the estimation for $E_Q|N^{(n,p)}_T |$.

Then, the Lebesgue dominated convergence theorem applied to the right-hand side of the previous inequality shows that it tends to zero as $n\rightarrow +\infty$.
On the other hand, from the fact that the Hellinger process is finite and also from the inequality $(\sqrt{Y(y)}-1)^2\geq Y(y)/6$ satisfied on $ A^ c$
we get
$$\int_{A^ c}Y(y) d\nu < +\infty$$
This result with the previous convergence  proves the convergence of $E_Q|N^{(n,p)}_T |$ to zero as $n\rightarrow +\infty$. Theorem 4 is proved.
$\Box$. 

\section{Utility maximising strategies}
We are now going to join the decomposition of the previous section with Theorem 3.1 of \cite{GR} in order to get an explicit expression of the optimal strategy.
                                 
\begin{theo} \label{strat}                                     
Let $u$ be a $\mathcal{C}^3(]\underline{x},+\infty[)$ utility function and $f$ its convex conjugate. Assume there exists an $f$-minimal martingale measure $Q^*$ which preserves the Levy property and such that the integrability conditions (\ref{integcd}) are satisfied.\\ 
Then for any fixed initial capital  $x> \underline{x}$, there exists a predictable process $\hat{\phi}$ and a sequence of stopping times $(\tau_n)_{n\geq 1}$ such that $(\hat{\phi}_{.\wedge \tau_n})_{n\geq 1}$ is asymptotically $u$-optimal. In addition, $\hat{\phi}$ defines a $u$-optimal strategy as soon as $\underline{x}>-\infty$. \\
Furthermore, if $c\neq 0$, we have 
$$\hat{\phi}^{(i)}_s=-\frac{\lambda \beta^{(i)} Z_{s-}}{S^{(i)}_{s-}}\,\xi _s(\lambda Z_{s-})$$
where $\beta = ^{\top}(\beta ^{(1)},\cdots, \beta ^{(d)})$ is the first Girsanov parameter, the process $\xi_s(\cdot)$ is defined by (\ref{xi}) and $\lambda$ is a unique solution to the equation $E_{Q^*}(-f'(\lambda Z_T))=x$.\\
If $c=0$ and $\stackrel{\circ}{supp}(\nu)\neq \emptyset$ and it contains zero, then $f''(x)= ax^{\gamma}$ with $a>0$ and $\gamma \in \mathbb R$, and
$$\hat{\phi}^{(i)}_s=-\frac{\lambda \gamma^{(i)} Z_{s-}}{S^{(i)}_{s-}}\,\xi _s(\lambda Z_{s-})$$
where again  $\lambda$ is a unique solution to the equation $E_{Q^*}(-f'(\lambda Z_T))=x$ and  the constants $\gamma^{(i)}$ are related with the second Girsanov parameter $Y$ by the formula:
\begin{equation}\label{gamma}
\gamma^{(i)} =\exp (-y_{0,i})\,Y(y_0)^{\gamma}\,\frac{\partial}{\partial y_i} Y(y_0)
\end{equation}
where $y_0$ is chosen arbitrarily in $\stackrel{\circ}{supp}(\nu)$.\\
\end{theo}

\vskip 0.2cm 
\noindent \it Proof of Theorem \ref{strat} \rm The first part of the Theorem is a slight adaptation of \cite {Ka}. We do however recall the proof for the reader's ease and because of some changes due to the use of asymptotically optimal strategies.\\ 
 As $f'$ is strictly increasing, continuous and due to (\ref{integcd}), the function $\lambda\mapsto E_{Q^*}[f'(\lambda Z_T)]$ is also strictly increasing and continuous. Furthermore, since $f'= -(u')^{-1}$, we have $\lim_{\lambda\to 0}E_{Q^*}[f'(\lambda Z_T)]=-\infty$ and $\lim_{\lambda\to+\infty}E_{Q^*}[f'(\lambda Z_T)]=\underline{x}$. Hence, for all $x>\underline{x}$, there exists a unique $\lambda>0$ such that $-E_{Q^*}[f'(\lambda Z_T)]=x$. As $Q^*$ is minimal for the function $x\mapsto f(\lambda x)$, it follows from Theorem 3.1 of \cite{GR}, that there exists a predictable process $\hat{\phi}$ such that 
\begin{equation}\label{eqphi1}
-f'(\lambda Z_T)=x+(\hat{\phi}\cdot S)_T
\end{equation}
and furthermore $\hat{\phi}\cdot S$ defines a $Q^*$-martingale. By definition of the convex conjugate, we have 
$$u(x+(\hat{\phi}\cdot S)_T)=f(\lambda Z_T) +Z_T\,f'(\lambda Z_T)$$ and, hence,
$$E_P[|u(x+(\hat{\phi}\cdot S)_T|]\leq E_P|f(\lambda Z_T)| + E_P[Z_T|f'(\lambda Z_T)|]<+\infty.$$
If now $\phi$ denotes any admissible strategy, we have by definition of $f$, 
$$\begin{aligned}
u(x+(\phi\cdot S)_T)\leq &(x+(\phi\cdot S)_T)\lambda Z_T+f(\lambda Z_T)
\\\leq &(x+(\phi\cdot S)_T)\lambda Z_T+u(x+(\hat{\phi}\cdot S)_T)-\lambda Z_Tf'(\lambda Z_T)
\end{aligned}$$
and taking expectation, we obtain
$$E_P[u(x+(\phi\cdot S)_T)]\leq E_P[u(x+(\hat{\phi}\cdot S)_T]+\lambda E_{Q^*}[(\phi\cdot S)_T].$$
Now, under $Q^*$, $(\phi\cdot S)$ is a local martingale which is bounded from below, hence is a super-martingale, so that $E_{Q^*}[(\phi\cdot S)_T]\leq 0$. Therefore, 
$$E_P[u(x+(\phi\cdot S)_T)]\leq E_P[u(x+(\hat{\phi}\cdot S)_T)].$$
Furthermore, if $\underline{x}>-\infty$, we note that $(\hat{\phi}\cdot S)_T\geq \underline{x}-x$, so that $\hat{\phi}$ defines an admissible strategy, and hence is a $u$-optimal strategy. 
\par Now, using the concavity of the function $u$, we have that for all $n\geq 1$, 
$$u(x+(\hat{\phi}\cdot S)_{T\wedge \tau_n})\leq u(x+(\hat{\phi}\cdot S)_T)+u'(x+(\hat{\phi}\cdot S)_T)((\hat{\phi}\cdot S)_T-(\hat{\phi}\cdot S)_{T\wedge \tau_n}).$$
Taking absolute values, and using the fact that $u'=(-f')^{-1}$, we then have 
$$|u(x+(\hat{\phi}\cdot S)_{T\wedge \tau_n})|\leq |u(x+(\hat{\phi}\cdot S)_{T})|+Z_T|(\hat{\phi}\cdot S)_T-(\hat{\phi}\cdot S)_{T\wedge \tau_n}|.$$
As $(\hat{\phi}\cdot S)$ is a $Q^*$-uniformly integrable martingale, the family $((\hat{\phi}\cdot S)_{T\wedge \tau_n})_{n\geq 1}$ is uniformly integrable with respect to $P$. Hence, in particular, 
$$\lim_{n\to+\infty}E_P[u(x+(\hat{\phi}{\bf 1}_{[\![0, \tau_n]\!]}\cdot S)_T)]=E_P[u(x+\hat{\phi}\cdot S)_T]$$
Finally, $(\hat{\phi}{\bf 1}_{[\![0, \tau_n]\!]})$ is asymptotically $u$-optimal. 

\vskip 0.2cm We now want to obtain a more explicit expression for $\hat{\phi}$. First of all, we note that relation (\ref{eqphi1}) may be rewritten as 
$$-E_{Q^*}[f'(\lambda Z_T)|\mathcal{F}_t]=x+\sum_{i=1}^d \int_0^t \hat{\phi}^{(i)}_sS^{(i)}_{s-}dX^{(c)}_s+\int_0^t \int_{\mathbb{R}^{d*}}\hat{\phi}^{(i)}_sS_{s-}^{(i)}(e^{y_i} - 1)(\mu^{X}-\nu^{X})(ds, dy)$$
We can then identify this decomposition with that obtained in Theorem \ref{thdec}. If $c\neq 0$, we identify the continuous components and obtain that  $Q^*$-a.s, for all $t\leq T$, 
$$\sum_{i=1}^d \lambda\,\beta^{(i)}\int_0^t  \xi _s (\lambda Z_{s-})\,Z_{s-}\,dX^{(c),i}_s = - \sum_{i=1}^d \int_0^t\hat{\phi}^{(i)}_sS^{(i)}_{s-}dX^{(c),i}_s$$
Taking quadratic variation of the difference of the right and left-hand sides in the previous equality, we obtain that $Q^*$-a.s. for all $s\leq T$  
$$^{\top}[\lambda \beta\xi_s(\lambda Z_{s-})Z_{s-} + \phi_{s} \,S_{s-}]\,c\,[\lambda \beta\xi_s(\lambda Z_{s-})Z_{s-} + \phi_{s}\,S_{s-}] = 0$$
where by convention $\phi_{s}\,S_{s-}=(\phi^{(i)}_{s}\,S^{(i)}_{s-})_{1\leq i \leq d}$. 
Therefore, as $c$ is a symmetric positive matrix, we have 
$$\hat{\phi}_s S_{s-}= - \lambda \beta\xi_s(\lambda Z_{s-})Z_{s-} +V_s$$
where $V_s$ belongs to the kernel of $c$. We may now write
$$-E_{Q^*}[f'(\lambda Z_T)|\mathcal{F}_t] = x- \sum_{i=1}^d \lambda \beta^{(i)}\int_0^t\xi_s(\lambda Z_{s-})\frac{dS^{(i)}_s}{S_{s-}^{(i)}}+\sum_{i=1}^d \int_0^t V_s^{(i)}d X^{(c),i}_t$$
As for all $s\geq 0$,  $cV_s=0$, we must have $<\sum_{i=1}^d \int_0^{\cdot} V_s^{(i)}d X^{(c),i}>_s=0$, and so $Q^*$-a.s.
$$-E_{Q^*}[f'(\lambda Z_T)|\mathcal{F}_t] = x - \sum_{i=1}^d \lambda \beta^{(i)}\int_0^t \xi_s(\lambda Z_{s-})\frac{dS^{(i)}_s}{S_{s-}^{(i)}}$$
It then follows from the first part of the proof that the process $\hat{\phi}$ defined in (\ref{eqphi1}) defines an (asymptotically-) optimal strategy. 

\par If we now assume that $c=0$, we identify the discontinuous components and obtain that $Q^*$-a.s, for all $s\leq T$ and for  a.e. $y\in supp(\nu)$,
\begin{equation}\label{equal}
\sum_{i=1}^d \hat{\phi}^{(i)}_s S_{s-}^{(i)}(e^{y_i} -1)= - H_s( \lambda Z_{s-}, y)
\end{equation}
In addition, since  $\stackrel{\circ}{supp}(\nu)\neq \emptyset$ we obtain from Theorem 3 of \cite{CV} that
$$f'(xY(y))-f'(x)=\Phi(x)\sum_{i=1}^d \alpha^{(i)} (e^{y_i}-1)$$
where 
$$\Phi(x) = xf''(xY(y_0)),\,\,\, \alpha ^{(i)}= \exp(-y_{0,i}) \frac{\partial}{\partial y_i} Y(y_0)$$
with any $y_0\in \stackrel{\circ}{supp}(\nu)$. Again from  Theorem 5 of \cite{CV}, $f''(x) = a x^{\gamma}$  and this implies  after the derivation of (\ref{equal}) with respect to $y_i$, the formula for optimal strategy.
$\Box$  
\vskip 0.2cm \noindent We finally give a unified expression of optimal strategies for all utility functions associated with common f-divergence functions.

\begin{prop}
Consider a Levy process $X$ with characteristics $(b,c,\nu)$ and let $f$ be a function such that $f''(x)=ax^{\gamma}$, where $a>0$ and $\gamma\in\mathbb{R}$. Let $u_f$ be its concave conjugate. Assume there exists $\beta\in\mathbb{R}^d$ and a measurable function $Y:\mathbb{R}^ \setminus \{0\}\rightarrow \mathbb{R}^{+,d}$ such that 
\begin{equation}  \label{Y}
  Y(y)=(f')^{-1}(f'(1)+\sum_{i=1}^d \beta ^{(i)}(e^{y_i}-1))
\end{equation}
 and such that the following properties hold:
\begin{equation}\label{cdsec1}
Y(y)> 0 \,\,\,\nu-a.e.,
\end{equation}
\begin{equation}\label{cdsec2}
\sum_{i=1}^d \int_{|y|\geq 1}(e^{y_i}-1)Y(y)\nu(dy)<+\infty.
\end{equation}
\begin{equation} \label{cdsec3}
b+\frac{1}{2}diag(c)+c\beta+\int_{\mathbb{R}^{d}}((e^y-1)Y(y)-h(y))\nu(dy)=0.
\end{equation}
Then there exists a sequence of asymptotically optimal strategies $(\hat{\phi}_{.\wedge \tau_n})_{n\geq 1}$ whose coordinates are given by 
$$\hat{\phi}^{(i)}_s=\alpha_{\gamma+1}(x)\frac{\beta^{(i)}}{E_{Q^*}[Z^{\gamma+1}_{s}]}\frac{Z^{\gamma+1}_{s-}}{S^{(i)}_{s-}},$$
where $Z$ is the density process of the change of measure from P into the $f$-minimal equivalent martingale measure $Q^*$ and 
\begin{equation}\label{alpha}
\alpha_{\gamma+1}(x)= \frac{\gamma+1}{a}(x+f'(1))-1.
\end{equation}
In addition, $\hat{\phi}$ is optimal as soon as $\gamma\neq -1$. 
\end{prop}
\vskip 0.2cm \it Proof \rm We recall from \cite{CV} that under the assumptions (\ref{cdsec1}), (\ref{cdsec2}) and (\ref{cdsec3}), the Levy model has an $f$-minimal martingale measure which preserves the Levy property and whose Girsanov parameters are $(\beta,Y)$  if $c\neq 0$, and $(0,Y)$ if $c=0$. 
Let $\lambda>0$ be such that $E_{Q^*}[f'(\lambda Z_T)]=-x$. It is easy to see that in this case, the decomposition of Theorem \ref{thdec} can be written 
$$-E_{Q^*}[f'(\lambda Z_T)|\mathcal{F}_t]=x-a\lambda^{\gamma+1}\sum_{i=1}^d \beta^{(i)}\int_0^T Z^{\gamma+1}_{s-}\,\frac{E_{Q^*}[Z^{\gamma+1}_{T-s}]}{E_{Q^*}[Z^{\gamma+1}_{T}]}\frac{dS^{(i)}_s}{S^{(i)}_{s-}}$$
As $Q^*$ preserves the Levy property, we have $E_{Q^*}[Z^{\gamma+1}_{T-s}]\,E_{Q^*}[Z^{\gamma+1}_s]=E_{Q^*}[Z^{\gamma+1}_T]$, so that for $t=T$, this may be rewritten $$-f'(\lambda Z_T)=x+\alpha_{\gamma+1}(x) \sum_{i=1}^d \beta^{(i)}\int_0^T \frac{Z^{\gamma+1}_{s-}}{E_{Q^*}[Z^{\gamma+1}_{s}]}\frac{dS^{(i)}_s}{S^{(i)}_{s-}}$$
It then follows from the proof of Theorem \ref{strat} that $\hat{\phi}$ defines a sequence of asymptotically optimal strategy. $\Box$

\noindent \section{Acknowledgements} \rm This work was supported in part by ECOS project M07M01  and ANR-09-BLAN-0084-01 of Auto-similarity of  Department of Mathematics of Angers University.


\end{document}